\documentclass[11pt]{amsart}
\usepackage{amssymb}
\usepackage{amsmath}
\usepackage{amsfonts}
\newtheorem{theorem}{Theorem}
\theoremstyle{plain}
 \newcommand{\fin}{\hfill  $\square$}

\newtheorem {lemma}[theorem]{Lemma}
\newtheorem {prop}[theorem]{Proposition}

\newtheorem {remarks}[theorem]{Remarques}

\newtheorem {cor}[theorem]{Corollary}
\begin{document}

\title[ ]{The Multivariable moment problems and recursive relations }

\author{ K. Idrissi }
\address{  K. Idrissi, Mohamed V University. Rabat Morocco }
\email{kaissar.idrissi@gmail.com }

\author{E. H. Zerouali}
\address{E. H. Zerouali, Mohamed V University. Rabat Morocco}
\email{Zerouali@fsr.ac.ma}
\date{}

\begin{abstract}
  Let $\beta \equiv \{ \beta_\mathbf{i} \}_{\mathbf{i} \in \mathbb{Z}_+^d}$ be a $d$-dimensional multisequence.
   Curto and Fialkow, have shown that  if the infinite moment matrix $M(\beta)$ is finite-rank positive semidefinite, then $\beta$ has a unique representing measure, which is $rank M(\beta)$-atomic. Further, let $\beta^{(2n)} \equiv \{ \beta_\mathbf{i} \}_{\mathbf{i} \in \mathbb{Z}_+^d, \mid \mathbf{i} \mid \leq 2n}$ be a given truncated multisequence, with associated moment matrix $M(n)$ and $rank M(n)=r$, then $\beta^{(2n)}$ has an $r$-atomic representing measure $\mu$ supported in the semi-algebraic set $K=\{ (t_1, \ldots, t_d) \in \mathbb{R}^d : q_j(t_1, \ldots, t_d) \geq 0, 1\leq j\leq m \}$, where $q_j \in \mathbb{R}[t_1, \ldots, t_d]$, if $M(n)$ admits a positive rank-preserving extension $M(n+1)$ and the localizing matrices $M_{q_j}(n +[\frac{\deg q_j +1}{2}])$ are positive semidefinite; moreover, $\mu$ has precisely $rank M(n) - rank M_{q_j}(n +[\frac{\deg q_j +1}{2}])$ atoms in $\mathcal{Z}(q_j) \equiv \{ t\in \mathbb{R}^d: q_j(t)=0 \}$. In this paper, we show that every truncated moment sequence $\beta^{(2n)}$ is a subsequence of an infinite recursively generated multisequence, we investigate such sequences  to give an alternative proof of  Curto-Fialkow's results and also to obtain a new interesting results.
\end{abstract}

\maketitle

\section{ Introduction}

 Given a real multisequence $\{\beta_\mathbf{i}\}_{ \mathbf{i} \in \mathbb{Z}_+^d, \mid \mathbf{i} \mid \leq 2n}$ and a nonempty subset $K$ of $\mathbb{R}^d$.
The truncated $K$-moment problem for $\{\beta_\mathbf{i}\}_{ \mathbf{i} \in \mathbb{Z}_+^d, \mid \mathbf{i} \mid \leq 2n}$ entails finding a positive Borel measure $\mu$,  supported in $K$  such that

 \begin{equation}\label{i.1}
    \beta_\mathbf{i}= \int \mathbf{x}^{\mathbf{i}} d\mu   \phantom{text} ( \mathbf{i} \in \mathbb{Z}_+^d, \mid \mathbf{i} \mid \leq 2n ).
 \end{equation}
Let $\mathbf{i}$ be a multi-index $(i_1, i_2, \ldots, i_d)$ of positive integers. We will write
$ \mid \mathbf{i} \mid = i_1+ \ldots+ i_d$ and $\mathbf{x^i}= x_1^{i_1} \ldots  x_d^{i_d}$ whenever $\mathbf{x}\equiv (x_1, \ldots, x_d)$ is a $d$-tuple of real numbers.

 A solution $\{\beta_\mathbf{i}\}_{ \mathbf{i} \in \mathbb{Z}_+^d, \mid \mathbf{i} \mid \leq 2n}$ of \ref{i.1} is called a truncated moment sequence and $\mu$ is said to be a $K$-representing measure for $\{\beta_\mathbf{i}\}_{ \mathbf{i} \in \mathbb{Z}_+^d, \mid \mathbf{i} \mid \leq 2n}$. The full $K$- moment problem prescribes moments of all orders. More precisely, an infinite multisequence $\{\beta_\mathbf{i}\}_{ \mathbf{i} \in \mathbb{Z}_+^d}$ is given and we aim to find a positive Borel measure $\mu$ supported in $K$ such that

 \begin{equation}\label{i.2}
    \beta_\mathbf{i}= \int \mathbf{x}^{\mathbf{i}} d\mu \phantom{text} \text{for all } \mathbf{i} \in \mathbb{Z}_+^d;
 \end{equation}
 If $K= \mathbb{R}^d$; \eqref{i.2} is often, and we would, called the moment problem for $\{\beta_\mathbf{i}\}_{ \mathbf{i} \in \mathbb{Z}_+^d}$.

 In view of its fundamental importance in various field of mathematics and applied science, the results of $K$-moment problem have interesting application; for instant, the Curto-fialkow's results (\cite[Theorems 4.7]{CuFi1} and \cite[Theorems 1.6]{CuFi4}) have been crucial in the Lasserre's method for minimizing a polynomial over a semialgebraic set (which is NP-hard in general), see for instance \cite{Las1, Las2, Las3}. This wide area of applications   motivated Curto an Fialkow to give a generalization in several variable.

 \begin{theorem}\label{th1}\cite{CuFi10, CuFi4}
  Let $\beta \equiv \{\beta_\mathbf{i}\}_{ \mathbf{i} \in \mathbb{Z}_+^d}$  be a real multisequence and let $M(\beta)$ be its associated moment matrix. If  $M(\beta) \geq 0$ and $M(\beta)$ have  finite $rank$, then $\beta$ has a unique $rank M(\beta)$-atomic representing measure.
 \end{theorem}

  The next theorem characterizes  the existence of a finitely atomic $K$-representing measure, in the case when $K$ is a semi-algebraic set, that is,
 $K= K_{\mathcal{Q}}:= \{ (t_1, \ldots, t_d) \in \mathbb{R}^d : q_i(t_1, \ldots, t_d) \geq 0, 1\leq i\leq m \}$ where $\mathcal{Q}\equiv \{ q_i \}_{i=1}^m \subset \mathbb{R}[t_1, \ldots, t_d]$.

 \begin{theorem}\label{th2}\cite[Theorem 2.9]{CuFi10}
  An $d$-dimensional real sequence $\beta^{(2n)} \equiv \{\beta_\mathbf{i}\}_{ \mathbf{i} \in \mathbb{Z}_+^d, \mid \mathbf{i} \mid \leq 2n}$ admits a $rank M(n)$-atomic representing measure supported in $K_{\mathcal{Q}}$ if and only if $M(n) \geq 0$ and $M(n)$ admits a flat extension $M(n+1)$ such that $M_{q_i}(n+[\frac{\deg q_i +1}{2}]) \geq 0$ $(1\leq i\leq m)$.
  In this case, $M(n+1)$ admits a unique representing measure $\mu$, which is a $rank M(n)$-atomic (minimal) $K_{\mathcal{Q}}$-representing measure for $\beta$; moreover, $\mu$ has precisely $rank M(n) - rank M_{q_i}(n + [\frac{\deg q_i +1}{2}])$ atoms in $\mathcal{Z}(q_i) \equiv \{ t\in \mathbb{R}^d: q_i(t)=0 \}$, $1\leq i\leq m$.
 \end{theorem}

   An  alternative proof was provided  by M. Laurent in \cite{Lau}, based on Hilbert's Nullstellensatz instead of the functional analytic tools used in the original proof of Curto and Fialkow.

   The main purpose of this paper is to use the multi-indexed recursively sequences (and a  Binet formula) not only in the aim to obtain a new and short   proof of Theorem \ref{th1} and Theorem \ref{th2}, but also to give a new approach in solving  the K-moment problem.

\section{ On the recursively generated multisequences and moment matrix}
  In this section, we state two necessary conditions for the existence of a finitely atomic $K$-representing measure.   We will prove that these conditions are, also, sufficient for the existence of a, unique  minimal, representing measure.

  \subsection{ Two necessary conditions}

   The first necessary condition is established by studying the matrix positivity. To this end, we  introduce some notation and definitions.

Let $\beta \equiv \{\beta_\mathbf{i}\}_{ \mathbf{i} \in \mathbb{Z}_+^d, \mid \mathbf{i} \mid \leq 2n}$   be  a given moment multisequence associated with a representing measure $\mu \ge 0$.

The matrix
   $M(n)(\beta) \equiv M(n) \equiv M(\beta) := ( \beta_{\mathbf{i}+\mathbf{j}})$ is known in Curto-Fialkow's terminology as moment matrix. The  columns and rows, of $M(n)$, are  labeled  by the lexicographic ordering of the canonical basis of the real vector space $\mathbb{R}_n[x_1, \ldots, x_d]$ of real-valued polynomials of degree at most $n$.  For example, in  the case where  $d=2$ we write
       $$ 1; X_1; X_2; X_1^2; X_1 X_2; X_2^2;  \cdots X_1^n; X_1^{n-1} X_2;\cdots, X_2^n.$$
  Clearly,  the entry of $M(n)$ in row $\mathbf{X}^{\mathbf{i}}$ and column $\mathbf{X}^{\mathbf{j}}$  is $M(n)_{\mathbf{ij}}= \beta_{\mathbf{i}+\mathbf{j}}$. Furthermore, let $P \in \mathbb{R}[x_1, \ldots, x_d]$ with coefficient vector $\{ P_\gamma \}$ and let $P*\beta$ denote the   vector in $\mathbb{R}^{\mathbb{Z}_+^d}$ whose $\alpha$-th entry is $(P*\beta)_\alpha := \sum\limits_{\gamma}P_\gamma \beta_{\gamma +\alpha}$. The moment matrix
   $M(n)(P*\beta) = M_P (n+ [ \frac{1+\deg P}{2} ])$ is called the localizing matrix with respect to $\beta$ and $P$.

   Let  $\mu$ be a positive Borel measure on $\mathbb{R}^d$ such that
   \begin{equation*}
   \beta_{\mathbf{i}}= \int \mathbf{x^i} d\mu, \phantom{ text } \mathbf{i} \in \mathbb{Z}_+^d.
   \end{equation*}
   Then given $p(\mathbf{x})= \sum\limits_{\mid \mathbf{i}\mid \leq n} a_\mathbf{i} \mathbf{x^i} \in \mathbb{R}_n[x_1, \ldots, x_d]$, we have
   \begin{equation*}
   0 \leq \int p^2(\mathbf{x}) d\mu = \int \sum\limits_{\mathbf{i}, \mathbf{j}} a_{\mathbf{i}} a_{\mathbf{j}} \mathbf{x}^{\mathbf{i}+\mathbf{j}} d\mu
      = \sum\limits_{\mathbf{i}, \mathbf{j}} a_{\mathbf{i}} a_{\mathbf{j}} \int \mathbf{x}^{\mathbf{i}+\mathbf{j}} d\mu
      = \sum\limits_{\mathbf{i}, \mathbf{j}} a_{\mathbf{i}} a_{\mathbf{j}} \beta_{\mathbf{i}+\mathbf{j}},
   \end{equation*}

   hence $M(n)(\beta)$ is positive semi-definite.

   Moreover, if $\mu$ is supported on the closed semialgebraic set $F= \{ \mathbf{x} \in \mathbb{R}^d \mid h_1(\mathbf{x})\geq 0, \ldots, h_m(\mathbf{x})\geq 0 \}$, where $h_j \in \mathbb{R}[x_1, \ldots, x_d]$, then

   \begin{equation*}
   \begin{aligned}
0\leq &\int p^2(\mathbf{x}) h_l(\mathbf{x}) d\mu \\
 &= \sum\limits_{\mathbf{i},\mathbf{j}} a_{\mathbf{i}} a_{\mathbf{j}} \int \mathbf{x}^{\mathbf{i}+\mathbf{j}} h_l(\mathbf{x}) d\mu \\
 &= \sum\limits_{\mathbf{i},\mathbf{j}} a_{\mathbf{i}} a_{\mathbf{j}} \int \sum\limits_\alpha (h_l)_{\mathbf{\alpha}} \mathbf{x}^{\mathbf{i}+\mathbf{j}+\mathbf{\alpha}} d\mu \\
 &= \sum\limits_{\mathbf{i},\mathbf{j}} a_{\mathbf{i}} a_{\mathbf{j}} (\sum\limits_\alpha (h_l)_{\mathbf{\alpha}} \beta_{\mathbf{i}+\mathbf{j}+\mathbf{\alpha}})\\
 &= \sum\limits_{\mathbf{i},\mathbf{j}} a_{\mathbf{i}} a_{\mathbf{j}} (h_l*\beta)_{\mathbf{i}+\mathbf{j}}
 ,  \text{ for all } l=0, \ldots,m.
   \end{aligned}
   \end{equation*}
  It follows that, if $\beta \equiv \{\beta_\mathbf{i}\}_{ \mathbf{i} \in \mathbb{Z}_+^d, \mid \mathbf{i} \mid \leq 2n}$ admits a representing measure supported on $F$, then the matrices $M_{h_1}(n +[\frac{\deg h_1 +1}{2}]), \ldots, M_{h_m}(n +[\frac{\deg h_m +1}{2}])$ and $M(n)$ are semidefinite positive.

  The second necessary condition gives rise to recurrence relations. Indeed, suppose that $\beta \equiv \{\beta_\mathbf{i}\}_{ \mathbf{i} \in \mathbb{Z}_+^d, \mid \mathbf{i} \mid \leq 2n}$ is a moment sequence, a result of C. Bayer and J. Teichmann \cite{BT} states  that every finite moment sequence admits a finite atomic representing measure. Hence there exists a representing measure $\mu$ for $\beta$ such that
  $\text{supp} \mu \subseteq \times_{i=1}^{m_l} \{ \lambda_{l, 1}, \ldots, \lambda_{l, m_l}\}$. We
  write $p_l(x) =  \prod_{l=1}^d (x_l- \lambda_{l, i}) = x_l^{m_l} -a_1^{(l)} x_l^{m_l -1} -\ldots -a_{m_l}^{(l)}$, for all $l \in \{1, \ldots, d\}$, we have
  \begin{equation*}
  \begin{aligned}
  0 &= \int \mathbf{x^i} p_l(x) d\mu, \phantom{ texttttt} \text{for all } \mathbf{i} \in \mathbb{Z}_+^d, \\
    &= \beta_{m_l \epsilon_l +\mathbf{i}} -a_1^{(l)} \beta_{(m_l -1)\epsilon_l +\mathbf{i}} - \ldots - a_{m_l}^{(l)} \beta_{\mathbf{i}},
  \end{aligned}
  \end{equation*}
 where $\epsilon_l$ is the $d$-tuple   with $1$ in the $l$-th place and zero elsewhere. It follows  that every truncated moment sequence  can be regarded as the initial data  of an infinite moment  sequence verifying the following recurrence relations,
 \begin{equation}\label{i.3}
 \beta_{(m_l +1) \epsilon_l +\mathbf{i}} = a_0^{(l)} \beta_{m_l \epsilon_l +\mathbf{i}} + a_1^{(l)} \beta_{(m_l -1)\epsilon_l +\mathbf{i}}+ \ldots +a^{(l)}_{m_l}\beta_\mathbf{i},
 \end{equation}
for ever $\mathbf{i} \in \mathbb{Z}_+^d \text{ and } l \in \{1, \ldots, d\}$.

  It results that the recursiveness is inherent in the truncated moment problem. This is our main  motivation in the study of  sequences satisfying  \eqref{i.3} in  connection with the moment problem.

  \subsection{ Recursively generated multisequences }

   Let $\{ a_j^{(l)} \}_{1\leq l \leq d, 0\leq j\leq m_l}$ be some fixed real numbers and let $\beta \equiv \{ \beta_{\mathbf{i}}\}_{\mathbf{i} \in \mathbb{Z}^+_d}$ be a real multisequence defined by the following recurrence relations:

  \begin{equation}\label{2.1}
 \beta_{(m_l +1) \epsilon_l +\mathbf{i}} = a_0^{(l)} \beta_{m_l \epsilon_l +\mathbf{i}} + a_1^{(l)} \beta_{(m_l -1)\epsilon_l +\mathbf{i}}+ \ldots +a^{(l)}_{m_l}\beta_\mathbf{i} \phantom{aaa}  (1\leq l\leq d),
 \end{equation}
 where $\omega = \{ \beta_{\mathbf{i}}\}_{\mathbf{i} \in \times_{l=1}^d \{0, \ldots, m_l\} }$ are given initial conditions.\\
 In the sequel, we shall refer to such sequence as recursively generated multisequence associated with the characteristic polynomials $p_\beta  \equiv (p_1, \ldots, p_d)$, where $$p_l(x) = x^{m_l +1} -a_0^{(l)} x^{m_l} -a_1^{(l)} x^{m_l -1} -\ldots -a_l^{(l)} \in \mathbb{R}[x]\: \: \: \:  (1 \le l \le d).$$

  A recursively generated multisequence can be defined in various ways using different characteristic polynomials as is shown in the following example. Let $\{ \beta_{(n, m, v)} \}_{(n, m, v) \in \mathbb{Z}^3_+}$ with $\beta_{(n, m, v)} = 5^m a^n (2^v -1)$, where $a$ is a nonzero real number. Then    $p_\beta = (x-a, x^2 +ax -5x -5a, x^2 -x -2)$ and $p'_\beta = (x^2-a^2, x -5, x^3 -2x^2 -x +2)$ are both characteristic polynomials of $\beta$.

  Let $\mathcal{P}_\beta$ denote the set of characteristic polynomials associated with $\beta$.

  \begin{remarks}
  \begin{enumerate}
    \item For every  $p_\beta \equiv (p_1, \ldots, p_d) \in \mathcal{P}_\beta$ and for  every
 $Q_1(x), \ldots, Q_d(x) \in \mathbb{R}[X]$, we have $(p_1 Q_1, \ldots, p_d Q_d) \in \mathcal{P}_\beta$.
    \item  The characteristic polynomials $p_\beta$, together with the initial conditions, are said to define the sequence $\beta$.
  \end{enumerate}
  \end{remarks}

  For reason of simplicity, we identify a polynomial $p \equiv \sum\limits_{\mid\mathbf{i}\mid \leq n} a_{\mathbf{i}} \mathbf{x}^{\mathbf{i}}$ with its coefficient vector $p = (a_{\mathbf{i}})$ with respect to the basis of monomials of $\mathbb{R}_n[x_1, \ldots, x_n]$ in degree-lexicographic order. Clearly, for every polynomials  $p \equiv \sum\limits_\mathbf{i} a_{\mathbf{i}} \mathbf{x}^{\mathbf{i}},
   q \equiv \sum\limits_\mathbf{j} b_{\mathbf{j}} \mathbf{x}^{\mathbf{j}} \in \mathbb{R}_n[x_1, \ldots, x_n]$, we have
  $p^T M(\beta) q = \sum\limits_{\mathbf{i}, \mathbf{j}} a_{\mathbf{i}} b_{\mathbf{j}} \beta_{\mathbf{i} +\mathbf{j}}$.

    The next lemma is an immediate consequence of \eqref{2.1}.
  \begin{lemma}\label{l.0}
  Let $\beta \equiv \{ \beta_{\mathbf{i}}\}_{\mathbf{i} \in \mathbb{Z}^+_d}$ be a recursively generated multisequence and let $M(\beta)$ be its associated moment matrix. Then $(p_1, \ldots, p_d)$ is a characteristic polynomials of $\beta$ if and only if $M(\beta) p_l=0$ (for all $l = 0, 1, \ldots,d$).
  \end{lemma}

  Singly indexed  sequences $S \equiv \{s_k \}_{k\in\mathbb{Z}_+}$ verifying \eqref{2.1}, with $d=1$, are   known in literature as weighted generalized
 Fibonacci sequence \cite{BRZ, DRS}.\\

   \begin{theorem}\cite[Theorem 1]{DRS}(Binet formula)
    Let $S \equiv \{s_k \}_{k\in\mathbb{Z}_+}$ be a generalized Fibonacci sequence, associated with the characteristic polynomial
    $p(x) = \prod\limits_{i=0}^{n-1} (x -\lambda_i)^{k_i}$, then
 \begin{equation}\label{binet}
 s_k = \sum\limits_{i=0}^{n-1} \sum\limits_{j=0}^{k_i -1} c_{i,j} k^j \lambda_i^k \phantom{text} (c_{i, k_i} \neq 0),
 \end{equation}
 where the $c_{i, j}$ are determined by the initial condition $s_k = 0, 1, \ldots, \deg p -1.$
   \end{theorem}

  Let us observe that if the characteristic polynomial of $S \equiv \{s_k \}_{k\in\mathbb{Z}_+}$ has distinct roots, say $p(x) = \prod\limits_{i=0}^{n-1} (x -\lambda_i)$, then \eqref{binet} can be written as follows:
  \begin{equation}\label{binet.1}
 s_k = \sum\limits_{i=0}^{n-1} c_{i} \lambda_i^k \phantom{text} ,
 \end{equation}
 where the $c_i$ are determined by the initial condition $s_k = 0, 1, \ldots, n -1.$

  As observed in \cite[Proposition 2.1]{BRZ} in the singly indexed case, among all characteristic polynomials defining  $S$, there exists a unique monic characteristic polynomial $p_S$ of minimal degree, called the minimal characteristic polynomial, and which divides every characteristic polynomial. The next proposition gives a generalization of this result.

  \begin{prop}\label{p.1}
  For every recursively generated multisequence $\beta \equiv \{ \beta_{\mathbf{i}}\}_{\mathbf{i} \in \mathbb{Z}^+_d}$ given by \eqref{2.1}, there exists  unique monic  characteristic polynomials $p_\beta= (p_1^{(\beta)}, \ldots, p_d^{(\beta)}) \in \mathcal{P}_\beta$ with minimal degree. Moreover, for all $(Q_1, \ldots, Q_d) \in \mathcal{P}_\beta$, $Q_l$ is a multiple of $p_l$ whenever $l\in \{ 1, \ldots, d\}$.
  \end{prop}

  {\it Proof.}  For $l \in \{ 1, \ldots, d \}$ and   $\mathbf{I}_l = (i_1, \ldots, i_{l-1}, i_{l+1}, \ldots, i_d) \in \mathbb{Z}^{d-1}_+$, a fixed $(d-1)$-tuple, we have
  \begin{equation*}
  \beta_{(i_1, \ldots, i_{l-1},m_l +1 +i_l, i_{l+1}, \ldots, i_d)} = a_0^{(l)} \beta_{(i_1, \ldots, i_{l-1},m_l +i_l, i_{l+1}, \ldots, i_d)} +\ldots +a_{m_l}^{(l)} \beta_{(i_1, \ldots, i_l, \ldots, i_d)},
  \end{equation*} for every  $i_l \in \mathbb{Z}_+^d$.\\

 Hence $p_l(x)= x^{m_l+1} -a_0^{(l)} x^{m_l} -\ldots -a_{m_l}^{(l)}$ is a characteristic polynomial associated with the general Fibonacci sequence
 $\beta^{(l)}: i_l \rightarrow \beta_(i_1, \ldots, i_l, \ldots, i_d)$. Thus there exists a minimal characteristic polynomial $p_{\beta, \mathbf{I}_l}$ associated with $\beta^{(l)}$.  Now, for  $\mathbf{I}_l \in \mathbb{Z}_+^{d-1}$, $p_{\beta, \mathbf{I}_l}$ divides $p_l$ and this implies that the polynomial
 $p_l^{(\beta)}= \bigwedge\limits_{\mathbf{I}_l \in \mathbb{Z}_+^{d-1}} p_{\beta, \mathbf{I}_l}$, the smallest common multiple of all $p_{\beta, \mathbf{I}_l}$, provides a positive answer to the proposition. \fin

  In the   remainder of this paper, we will associate with every recursively generated multisequence $\beta$ with its minimal polynomial, that we denote  $p_\beta$.

  \begin{prop}\label{p.2}
  Let $\beta \equiv \{ \beta_{\mathbf{i}}\}_{\mathbf{i} \in \mathbb{Z}^+_d}$ be a recursively generated multisequence, associated with the characteristic polynomials $(p_1, \ldots, p_d)$. If $M( \beta) \geq 0$, then, for every $l \in \{ 1, \ldots, d\}$, the polynomial $p_l(x)$ has distinct roots.
  \end{prop}
To prove Proposition \ref{p.2}, we need the following two lemmas of independent interest
\begin{lemma}\label{p3.1} Under the notations above, for every  $f, g, h \in \mathbb{R}[x_1, \ldots, x_d]$, we have
  \begin{equation}
  f^T M(\beta) (gh) = (fg)^T M(\beta) h.
  \end{equation}
\end{lemma}

  {\it Proof.} Let  $f, g, h \in \mathbb{R}[x_1, \ldots, x_d]$ be polynomials. We write $f =\sum\limits_{\mathbf{i}} f_{\mathbf{i}} \mathbf{x}^{\mathbf{i}}$, $g =\sum\limits_{\mathbf{j}} g_{\mathbf{j}} \mathbf{x}^{\mathbf{j}}$ and
  $h =\sum\limits_{\mathbf{k}} h_{\mathbf{k}} \mathbf{x}^{\mathbf{k}}$. As the entry of the moment matrix corresponding to the column $\mathbf{x}^{\mathbf{i}}$ and the line $\mathbf{x}^{\mathbf{j}}$ is $\gamma_{\mathbf{i}+\mathbf{j}}$, we obtain
  $$\begin{array}{lll}
f^T M(\beta) (gh)
&  = & (\sum\limits_{\mathbf{i}} f_{\mathbf{i}} \mathbf{x}^{\mathbf{i}})^T M(\beta) (\sum\limits_{\mathbf{j}, \mathbf{k}} g_{\mathbf{j}} h_{\mathbf{k}}
\mathbf{x}^{\mathbf{j}+\mathbf{k}})\\
  &=& \sum\limits_{\mathbf{i}, \mathbf{j}, \mathbf{k}} f_{\mathbf{i}} g_{\mathbf{j}} h_{\mathbf{k}} \gamma_{ \mathbf{i}+\mathbf{j}+\mathbf{k} }\\
  &=& (fg)^T M(\beta) h. \end{array}$$

   \begin{lemma}\label{p.x}   For every  polynomial $p \in \mathbb{R}[x_1, \ldots, x_2]$ and any integer $n\geq 1$, we have
  \begin{equation}\label{p3.2}
  M(\beta) p^n =0  \Longrightarrow   M(\beta) p =0.
  \end{equation}
\end{lemma}

 {\it Proof.}    If $M(\beta) p^2 =0$, then $0 = 1^T M(\beta) p^2 = p^T M(\beta) p$, from Lemma \ref{p3.1}; since $M(\beta) \geq 0$,  we obtain $M(\beta) p = 0$ and hence \eqref{p3.2} holds for $n=2$. By induction, \eqref{p3.2} remains valid for any power of $2$. Now , if $M(\beta) p^n =0$ we choose $r$ in such a way that $r+k$ is a power of $2$, hence
  $$M(\beta) p^{n+r} =0.$$
Which gives $M(\beta) p =0$.\\

{ \it Proof of Proposition \ref{p.2}.}
  Let $p_l(x) = \prod\limits_{i=0}^{m_l -1} (x -\lambda_{l, i})^{n_{l, i}}$ and let $M_l = \max\limits_{i=0}^{m_l -1} n_{l, i}$. Applying Lemma \ref{l.0} we obtain $M(\beta) \prod\limits_{i=0}^{m_l -1} (x -\lambda_{l, i})^{n_{l, i}} =0$,  we derive that  $M(\beta) (\prod\limits_{i=0}^{m_l -1} (x -\lambda_{l, i}))^{M_l} =0$, and hence Lemma \ref{p.x} yields that  $M(\beta) \prod\limits_{i=0}^{m_l -1} (x -\lambda_{l, i}) =0$. It follows, by Lemma \ref{l.0}, that  $\prod\limits_{i=0}^{m_l -1} (x -\lambda_{l, i})$ is a characteristic polynomial, of $\beta$, dividing  $p_l(x)$. Since $p_l(x)$ is minimal, then  $p_l(x) = \prod\limits_{i=0}^{m_l -1} (x -\lambda_{l, i})$, as desired.

\section{ Main results}

 We present a characterization of moment sequences  involving the recursively generated multisequence and the moment matrix, that leads to   new  proofs  of Theorem \ref{th1} and  Theorem \ref{th2}.

 \begin{theorem}\label{11}
  Let $\beta \equiv \{ \beta_{\mathbf{i}} \}_{\mathbf{i} \in \mathbb{Z}_+^d}$ be a multisequence of real numbers and let $M(\beta)$ and  $M(n)$ be the moment matrices  associated with $\beta$  and  $\{ \beta_{\mathbf{i}} \}_{\mid \mathbf{i} \mid \leq 2n, \mathbf{i} \in \mathbb{Z}_+^d}$, respectively. The following are equivalent:
  \begin{enumerate}
    \item $M(\beta)$ is a finite-rank positive semidefinite matrix.
    \item $\beta$ is recursively generated, associated with the minimal characteristic polynomials $(p_1, \ldots, p_d)$ and $M(\tau) \geq 0$, with
          $\tau =\sum\limits_{i=1}^d (\deg p_i -1)$.
    \item $\beta$ has a unique representing measure, which is $rank  M(\beta)$-atomic.
  \end{enumerate}
  \end{theorem}

   For the proof of Theorem \ref{11} we will need the following auxiliary result which can be regarded as the Binet Formula for the recursively generated multisequences.

   \begin{lemma}\label{l.3}
  Let $\beta \equiv \{ \beta_{(i_1, \ldots, i_d)} \}_{(i_1, \ldots, i_d) \in \mathbb{Z}_+^d}$ be a real recursively generated multisequence associated with the characteristic polynomials $(p_1, \ldots, p_d)$, where  $p_j (x) = \prod\limits_{l=0}^{m_j-1} (x -\lambda_{j, l})$. Then, there exists  $c_{(l_1, \ldots, l_d)}$ real numbers,  determined by the initial conditions $\{ \beta_{(i_1, \ldots, i_d)} \}_{ i_j \in \{0, \ldots, m_j -1\} }$ such that
  \begin{equation}\label{rl1}
  \beta_{(i_1, \ldots, i_d)}= \sum\limits_{l_1 =0}^{m_1-1} \ldots \sum\limits_{l_d =0}^{m_d-1} c_{(l_1, \ldots, l_d)} \lambda_{1, l_1}^{i_1} \ldots \lambda_{d, l_d}^{i_d}.
  \end{equation}

   \end{lemma}

 {\it Proof.}
  For  $(i_2, \ldots, i_d) \in \mathbb{Z}_+^{d-1}$ given, the singly sequence $i_1 \rightarrow  \beta_{(i_1, \ldots, i_d)}$ is a general Fibonacci sequence associated with the characteristic polynomial $p_1 (x) = \prod\limits_{l=0}^{m_1 -1} (x -\lambda_{1, l})$. Then the Binet formula implies that

 \begin{equation}\label{d.1}
 \beta_{(i_1, \ldots, i_d)}= \sum\limits_{l_1 =0}^{m_1-1}  c_{l_1}^{(i_2, \ldots, i_d)} \lambda_{1, l_1}^{i_1},
 \end{equation}
 where the $c_{l_1}^{(i_2, \ldots, i_d)}$ are determined by the initial condition $\{ \beta_{(i_1, \ldots, i_d)} \}_{ 0\leq i_1 \leq m_1 -1}$.
 We claim that, for every integer $l_1 \geq 0$, the single sequence $i_2 \rightarrow  c_{l_1}^{(i_1, \ldots, i_d)}$ is a general Fibonacci sequence associated with the characteristic polynomial $p_2(x)$.

 Indeed, since  $i_2 \rightarrow  \beta_{(i_1, i_2, \ldots, i_d)}$ is a Fibonacci sequence associated with $p_2(x)$,  we obtain
 \begin{equation*}
 \beta_{(i_1, i_2 +m_2, i_3, \ldots, i_d)} -a_1^{(2)} \beta_{(i_1, i_2 +m_2 -1, i_3, \ldots, i_d)} -\ldots -a_{m_2}^{(2)} \beta_{(i_1, i_2, \ldots, i_d)} =0.
 \end{equation*}
 It follows  from \ref{d.1},
 \begin{equation*}
 \sum\limits_{l_1 =0}^{m_1-1} p_2( c_{l_1}^{(i_2, \ldots, i_d)} ) \lambda_{1, l_1}^{i_1} =0
 \end{equation*}
 where $ p_2( c_{l_1}^{(i_2, \ldots, i_d)} ) = c_{l_1}^{(i_2 + m_2, i_3, \ldots, i_d)} -a_1^{(2)} c_{l_1}^{(i_2 + m_2 -1, i_3, \ldots, i_d)} -\ldots
 -a_{m_2}^{(2)} c_{l_1}^{(i_2, i_3, \ldots, i_d)}$.
 For $i_1 \in \{0, 1, \ldots, m_1 -2\}$ and $m_1 -1$, we derive the following Vandermonde system
 \begin{equation*}\left\{\begin{aligned}
 &p_2( c_{0}^{(i_2, \ldots, i_d)} ) \lambda_{1, 0}^{0} &+\ldots &+p_2( c_{m_1 -1}^{(i_2, \ldots, i_d)} ) \lambda_{1, m_1 -1}^{0} &=0\\
 &\vdots &  &  \vdots &\vdots\\
 &p_2( c_{0}^{(i_2, \ldots, i_d)} ) \lambda_{1, 0}^{m_1 -1} &+\ldots &+p_2( c_{m_1 -1}^{(i_2, \ldots, i_d)} ) \lambda_{1, m_1 -1}^{m_1 -1} &=0.
 \end{aligned}\right.
 \end{equation*}
 Since $\{\lambda_{1, l_1}\}_{l_1=0}^{m_1 -1}$ are distinct, the unique solution is zero,
 \begin{equation*}
 p_2( c_{l_1}^{(i_2, \ldots, i_d)} ) =0,   \phantom{ du text } 0\leq l_1 \leq m_2 -1.
 \end{equation*}
 As the integer $i_2$ is arbitrary, we have , $i_2 \rightarrow c_{l_1}^{(i_2, \ldots, i_d)}$ is a general Fibonacci sequence associated with $p_2(x)$ for every   $0\leq l_1 \leq m_1 -1$.  Similarly,  one can show that the singly indexed sequence $i_j \rightarrow c_{l_1}^{(i_2, \ldots, i_j, \ldots, i_d)}$ is a general Fibonacci sequence associated with the characteristic polynomial $p_j(x)$. By applying the Binet formula to the sequence $i_2 \rightarrow c_{l_1}^{(i_2, \ldots, i_d)}$, we get
 $c_{l_1}^{(i_2, \ldots, i_j, \ldots, i_d)} = \sum\limits_{l_2=0}^{m_2 -1} c_{(l_1, l_2)}^{(i_3, \ldots, i_d)} \lambda_{2, l_2}^{i_2}$, where
 $c_{(l_1, l_2)}^{(i_3, \ldots, i_d)}$ are determined by the initial condition $\{ c_{l_1}^{(i_2, \ldots, i_d)} \}_{0\leq i_2 \leq m_2 -1}$. Hence
 $$\beta_{(i_1, \ldots, i_d)} = \sum\limits_{l_1=0}^{m_1-1} \lambda_{1, l_1}^{i_1} \sum\limits_{l_2=0}^{m_2-1} \lambda_{2, l_2}^{i_2} c_{(l_1, l_2)}^{(i_3, \ldots, i_d)} = \sum\limits_{l_1=0}^{m_1-1} \sum\limits_{l_2=0}^{m_2-1} \lambda_{1, l_1}^{i_1} \lambda_{2, l_2}^{i_2} c_{(l_1, l_2)}^{(i_3, \ldots, i_d)}.$$
 Now we will show that, for every $0\leq l_1 \leq m_1 -1$ and $0\leq l_2 \leq m_2 -1$, the sequence $i_3 \rightarrow c_{(l_1, l_2)}^{(i_3, \ldots, i_d)}$
 is a general Fibonacci sequence associated with $p_3(x)$. To this aim it suffices to remark that $i_3 \rightarrow c_{l_1}^{(i_2, \ldots, i_j, \ldots, i_d)}$
 is a general Fibonacci sequence associated with $p_3(x)$ and replace, in the above proof,  the sequence  $i_2 \rightarrow  \beta_{(i_1, i_2, \ldots, i_d)}$ by
 $i_3 \rightarrow c_{l_1}^{(i_2, \ldots, i_j, \ldots, i_d)}$. Therefore, we obtain
 \begin{equation*}
 \beta_{(i_1, \ldots, i_d)}= \sum\limits_{l_1 =0}^{m_1-1} \sum\limits_{l_2 =0}^{m_2-1} \sum\limits_{l_3 =0}^{m_3-1}  \lambda_{1, l_1}^{i_1} \lambda_{2, l_2}^{i_2} \lambda_{3, l_3}^{i_3}  c_{(l_1, l_2, l_3)}^{(i_4, \ldots, i_d)},
 \end{equation*}
 where $c_{(l_1, l_2, l_3)}^{(i_4, \ldots, i_d)}$ are determined by $\{ c_{(l_1, l_2)}^{(i_3, \ldots, i_d)} \}_{ 0\leq i_3 \leq m_3 -1}$.  By induction we get
  \begin{equation}\label{exp}
 \beta_{(i_1, \ldots, i_d)}= \sum\limits_{l_1 =0}^{m_1-1} \ldots \sum\limits_{l_d =0}^{m_d-1} c_{(l_1, \ldots, l_d)} \lambda_{1, l_1}^{i_1} \ldots \lambda_{d, l_d}^{i_d},
  \end{equation}
  where  $ c_{(l_1, \ldots, l_d)}$ are  real numbers.
  \fin

Using the multi index notations ${\bf i}=(i_1,\cdots,i_p)$, ${\bf l}=(l_1,\cdots,l_p)$ and 
${\bf \lambda }_{\bf l}=(\lambda_{1, \mathbf{l}}, \cdots,\lambda_{p, \mathbf{l}})$, The expression \eqref{exp} becomes
\begin{equation}\label{rl2}
\beta_{\bf i} =   \beta_{(i_1, \ldots, i_d)}= \sum\limits_{\mathbf{l} \in \mathfrak{I}} c_{\mathbf{l}} \lambda_{1, \mathbf{l}}^{i_1} \ldots \lambda_{d, \mathbf{l}}^{i_d}
= \sum\limits_{\mathbf{l} \in \mathfrak{I}} c_{\mathbf{l}} {\bf \lambda }^{\bf i}_\mathbf{l},
  \end{equation}

   where  $\mathfrak{I} \equiv \mathfrak{I}(\beta) := \{ \mathbf{l}\equiv (l_1, \ldots, l_d) \in \mathbb{Z}_+^d \mid c_{\mathbf{l}} \neq 0 \}$.

 {\it Proof of Theorem \ref{11}.} First we show the equivalence between the assertions $(1)$ and $(2)$.  Let $j\in \{ 1, \ldots, d\}$ and let
 $s_{j+1} =\inf \{ n\in \mathbb{N} \mid X_j^{n+1} = a_0^{(j)} X^n_j -a_1^{(j)} X^{n-1}_j -\ldots -a_n^{(j)}1 \text{ for some } a_0^{(j)}, \ldots, a_n^{(j)} \in \mathbb{R} \}$, where $X_j^i$ denotes the index of columns and rows of $M(\beta)$ as in the above section. For every $\mathbf{i} \in \mathbb{Z}_+^d$, we will have
 \begin{equation*}
  \beta_{\mathbf{i} +(s_j +1) \epsilon_j} = a_0^{(j)} \beta_{\mathbf{i} +s_j \epsilon_j} +\ldots +a_{s_j}^{(j)} \beta_{\mathbf{i}}.
 \end{equation*}
 Hence $\beta$ is recursively generated associated with the minimal characteristic polynomials $(p_1, \ldots, p_d)$,
 where $p_j(x) = x^{s_j +1} -a_0^{(j)} x^{s_j} -\ldots -a_{s_j}^{(j)} (j\in \{ 1, \ldots,d\})$. \\
Conversely, since $\beta$ is recursively generated multisequence,  every column $\mathbf{X^i}$ in $M(\beta)$, such that $i_l \geq \deg p_l$ for some $l \in \{1, \ldots, d\}$, is a linear combination of lower (power index) columns; more precisely,
 \begin{equation*}
 \mathbf{X}^{(l_1, \ldots, l_j, \ldots, _d)} = \sum\limits_{i=1}^{m_l} a_i^{(j)} \mathbf{X}^{(l_1, \ldots, l_j -i, \ldots, _d)}.
 \end{equation*}
  and then $M(\beta)$ has a finite rank. It remains to show that $M(\beta) \geq 0$. To this end,  construct the  matrix $W_s \in M_{\tau +s +1, m(\tau +s +1)}$, the algebra of $(\tau +s +1)\times m(\tau +s +1)$ real matrices, where $m(\tau +s +1)$ denote the number of columns (or rows) of $M(\tau +s +1)(\beta)$, such that the successive columns  of $W_s$   are defined by
  \begin{equation*}
 \mathbf{X}^{\sum\limits_{k=1}^d l_k \epsilon_k}= \sum\limits_{i=1}^{\deg p_j} a_i^{(j)} e_{((l_j -i)\epsilon_j +\sum\limits_{k \neq j} l_k \epsilon_k)},
 \end{equation*}
 where $l_1+\ldots+l_d = \tau +s +1$,  $l_j \geq \deg p_j$ and $\{ e_\mathbf{i}\}_{\mid \mathbf{i}\mid \leq \tau +s +1}$ denote the canonical basis of $\mathbb{R}^{m(\tau +s +1)}$, that is, $e_\mathbf{i}$ is the vector with 1 in the $\mathbf{X}^\mathbf{i}$ entry and 0 all other positions.
  Remark that if $\sum\limits_{k=1}^d n_k \geq \tau +1$, then there exists $j \in \{1, \ldots, d\}$ such that $n_j \geq p_j$. Thus it follows, from $\ref{2.1}$, that
  \begin{equation*}
  M(\tau +s +1)=  \begin{pmatrix} M(\tau +s) &B\\ B^* &C \end{pmatrix},
  \end{equation*}
  with $B=M(\tau +s)W_{\tau +s}$ and $C=B^*W$. Therefore, if $M(\tau +s)\geq 0$, then (by Smul'jan's Theorem \cite{Smu}) we get $M(\tau +s +1) \geq 0$.
  As $M(\tau )\geq 0$ then, by induction over $s\geq 0$, we conclude that $M(\beta) \equiv M(\infty)(\beta)\geq 0$.

 We show now  the equivalence between $(2)$ and $(3)$. From \ref{p.2},  the polynomial $p_l(x)$ has distinct roots (not necessary real) for all $l \in \{1, \ldots, d \}$. We put  $p_l(x) = \prod\limits_{i=0}^{m_l -1} (x -\lambda_{l, i})$. According to the relation  \eqref{rl2}, $\beta_{\bf i} \equiv \beta_{(i_1, \ldots, i_d)}$ can be expressed as follows
  \begin{equation*}
   \beta_{\bf i} =   \beta_{(i_1, \ldots, i_d)}= \sum\limits_{\mathbf{l} \in \mathfrak{I}} c_{\mathbf{l}} \lambda_{1, \mathbf{l}}^{i_1} \ldots \lambda_{d, \mathbf{l}}^{i_d}
= \sum\limits_{\mathbf{l} \in \mathfrak{I}} c_{\mathbf{l}} {\bf \lambda }^{\bf i}_\mathbf{l}.
   \end{equation*}
   Thus the measure $$\mu = \sum\limits_{\mathbf{l} \in \mathfrak{I}} c_{\mathbf{l}}  d\delta_{\lambda_{1, \mathbf{l}}} \ldots d\delta_{\lambda_{d, \mathbf{l}}}= \sum\limits_{l_1 =0}^{m_1-1} \ldots \sum\limits_{l_d =0}^{m_d-1} c_{(l_1, \ldots, l_d)} d\delta_{\lambda_{1, l_1} } \ldots d\delta_{\lambda_{d, l_d}},$$
   satisfies  $$ \beta_\mathbf{i}= \int \mathbf{x}^{\mathbf{i}} d\mu.$$To see that $\mu$ provides  a positive answer to the moment problem \eqref{i.2}, that is, $c_{\mathbf{l}} > 0$ and $\lambda_{1, \mathbf{l}}, \ldots, \lambda_{d, \mathbf{l}} \in \mathbb{R}$, whenever $\mathbf{l} \in \mathfrak{I}$. We consider the the following family of  interpolation polynomials at the atoms of the representing measure $\mu$, say 
  $supp \mu := \{ \mathbf{\lambda}_{\mathbf{l}_1}, \ldots, \mathbf{\lambda}_{\mathbf{l}_m}\} \subset \mathbb{R}^d$,
  \begin{equation*}\begin{aligned}
  L_{\mathbf{\lambda}_{\mathbf{l}_s}} (x_1, \ldots, x_d) &= L_{(\lambda_{1, \mathbf{l}_s} , \ldots, \lambda_{d, \mathbf{l}_s})} (x_1, \ldots, x_d)\\
   &= \prod\limits_{i=0}^d (
  \tiny{\prod\limits_{\begin{matrix} 0\leq j\leq m \\ j\neq s \end{matrix}}} \frac{x_i -\lambda_{i, \mathbf{l}_j}}{\lambda_{i, \mathbf{l}_s} -\lambda_{i, \mathbf{l}_j}} ),   \phantom{aaa} (s \in \{1, \ldots, m\}).
  \end{aligned}
  \end{equation*}

  Clearly  
 $$L_{\mathbf{\lambda}_{\mathbf{l}_s}} (x_1, \ldots, x_d)=\{ \begin{array}{l} 1 \: \mbox{if}\: (x_1, \ldots, x_d) = (\lambda_{1, \mathbf{l}_s} , \ldots, \lambda_{d, \mathbf{l}_s}), \\
0 \: \mbox{elsewhere.}\end{array}$$

 It follows that, for any $\mathbf{l} \in \mathfrak{I}$,
  \begin{equation*}\begin{aligned}
  c_{\mathbf{l}} &= \int \mid L_{\mathbf{\lambda}_\mathbf{l}} \mid^2 d\mu \\
                         &= L_{\mathbf{\lambda}_\mathbf{l}}^T M(\beta) L_{\mathbf{\lambda}_\mathbf{l}} \geq 0,
  \end{aligned}
  \end{equation*}
 and also, for any $j \in \{0, \ldots, d \}$,
  \begin{equation*}\begin{aligned}
  \lambda_{j, \mathbf{l}}  c_{\mathbf{l}} &= \int x_j \mid L_{\mathbf{\lambda}_\mathbf{l}} \mid^2 d\mu \\
 &= L_{\mathbf{\lambda}_\mathbf{l}}^T M_{x_j}(\beta) L_{\mathbf{\lambda}_\mathbf{l}} \in \mathbb{R},
  \end{aligned}
  \end{equation*}
  since the localizing matrix $M_{x_j}(\beta)$ (defined above) is a symmetric real matrix. 
  As $c_{\mathbf{l}} \neq 0$, because $\mathbf{l} \in \mathfrak{I}$, then $c_{\mathbf{l}} > 0$, and hence  $\lambda_{j, \mathbf{l}} \in \mathbb{R}$, as desired.

It remains to show that $\mu$ is a $rank M(\beta)$-atomic and is the unique representing measure of $\beta$.
To this aim, assume that  $M(\beta) \sum\limits_{\mathbf{l} \in \mathfrak{I}} a_{\mathbf{l}} \frac{1}{ \sqrt{c_{\mathbf{l}}} } L_{\mathbf{\lambda_l}} = 0$,
 where $\{ a_{\mathbf{l}}\}_{\mathbf{l} \in \mathfrak{I}}$ are real numbers (not all zero) and $\mathbf{\lambda_l} = (\lambda_{1, \mathbf{l}} , \ldots, \lambda_{d, \mathbf{l}})$. Since $\frac{1}{\sqrt{c_{\mathbf{i}}}} L_{\mathbf{\lambda_i}}^T M(\beta) \frac{1}{\sqrt{c_{\mathbf{j}}}} L_{\mathbf{\lambda_j}} =\delta_{\mathbf{i}, \mathbf{j}}$, the Kronecker delta, we obtain $$0 = ( \sum\limits_{\mathbf{i} \in \mathfrak{I}} a_{\mathbf{i}} \frac{1}{\sqrt{c_{\mathbf{i}}}} L_{\mathbf{\lambda_i}} )^T M(\beta) \sum\limits_{\mathbf{i} \in \mathfrak{I}} a_{\mathbf{i}} \frac{1}{\sqrt{c_{\mathbf{i}}}} L_{\mathbf{\lambda_i}} = \sum\limits_{\mathbf{i} \in \mathfrak{I}} a_{\mathbf{i}}^2,$$ a contradiction. Thus $\text{card }supp \mu \leq rank M(\beta)$.

On the other hands, from \eqref{rl2},
 $M(\beta) = \sum\limits_{\mathbf{l} \in \mathfrak{I}} c_{\mathbf{l}} \zeta_{\mathbf{\lambda_i}}^T \zeta_{\mathbf{\lambda_i}}$, where
 $\zeta_{\mathbf{\lambda_l}} := (\mathbf{\lambda_l}^\alpha)_{\alpha \in \mathbb{Z}^d} \in \mathbb{R}^{\mathbb{Z}^d_+}$, hence
 $$rank M(\beta) \leq \text{card }\mathfrak{I} = \text{card }supp \mu.$$
 Therefore $rank M(\beta) = \text{card }supp \mu$.
 To get uniqueness, let us suppose now that $\mu' := \sum\limits_{\mathbf{i} \in \mathfrak{I}'} {c'}_\mathbf{i} \delta_{\mathbf{\lambda_i}}$ is another representing measure for $\beta$; that is, $\int p d\mu = \int p d\mu'$, for every $p \in \mathbb{R}[x_1, \ldots, x_d]$.
 Let $\{ L_{\mathbf{\lambda_i}} \}_{\mathbf{i} \in \mathfrak{I} \cup \mathfrak{I}'} \subset \mathbb{R}[x_1, \ldots, x_d]$ be the interpolating polynomials at the points of $\mathfrak{I} \cup \mathfrak{I}'$. If $supp \mu \neq supp \mu'$, then there exists $\mathbf{j} \in \mathfrak{I'} \backslash \mathfrak{I}$. Thus $$ 0\neq {c'}_\mathbf{j} = \int L_{\mathbf{\lambda_j}} d\mu' = \int L_{\mathbf{\lambda_j}} d\mu = 0,$$
 a contradiction, hence $supp \mu = supp \mu'$. Also, we have
 $${c'}_\mathbf{i} = \int L_{\mathbf{\lambda_i}} d\mu' = \int L_{\mathbf{\lambda_i}} d\mu =c_\mathbf{i}, \mbox{
 whenever }  \mathbf{i} \in \mathfrak{I} \cup \mathfrak{I}'.$$ Therefore $\mu = \mu'$, as desired.
 The reverse implication follows  directly from Section $2.1$.
 \fin

 Let us recall \cite[Theorem 7.8]{CuFi1}: If $M(n)$ is positive semidefinite and admits a flat extension $M(n+1)$, then $M(n+1)$ admits unique successive flat moment extensions $M(n+2), M(n+3), \ldots, M(\infty) \equiv M(\beta)$. In addition with Theorem \ref{11}, we obtain the following corollary.

\begin{cor}\label{c.1}
 Let  $\beta^{(2n)} \equiv \{\beta_\mathbf{i}\}_{ \mathbf{i} \in \mathbb{Z}_+^d, \mid \mathbf{i} \mid \leq 2n}$ be a truncated multisequence and let $M(n)$ be its associated moment matrix. If $M(n)$ is positive semidefinite and admits a flat extension $M(n+1)$, then $M(n+1)$ has a unique representing measure $\mu$; such that  $\text{card }supp \mu = rank  M(n)$.
\end{cor}

 We give now a short proof of Theorem \ref{th2}.

  {\it Proof of Theorem \ref{th2}.} We have shown in Section 2 that  the positive  semidefiniteness of $M(n)$ and $M_{q_i}(n +[\frac{\deg q_i +1}{2}])$, for all $q_i \in K_{\mathcal{Q}}$, are necessary conditions for the existence of a representing measure $\mu$, for $\beta^{(2n)}$, supported in $K_{\mathcal{Q}}$. Also, $\mu$ is a representing measure for some recursively generated moment sequence $\beta \equiv \{\beta_\mathbf{i}\}_{ \mathbf{i} \in \mathbb{Z}_+^d}$ (notice in passing  that, $M(\beta) \geq 0$ and $rank M(\beta) < \infty$). Therefore,  from  Theorem \ref{11}, $rank M(\beta) = supp \mu$ $(=  rank M(n))$ and hence $M(n)$ admits a rank-preserving extension.

We  prove the reverse inclusion. As $M(n)\geq 0$ and $M(n)$ admits a flat extension $M(n+1)$, from Corollary \ref{c.1}, $M(n+1)$ admits a unique $rank  M(n)$-atomic representing measure, write
  \begin{equation}\label{mu1}
  \mu = \sum\limits_{\mathbf{l} \in \mathfrak{I}} c_{\mathbf{l}} d\delta_{\mathbf{\lambda}_{1, \mathbf{l}}} \ldots d\delta_{\mathbf{\lambda}_{d, \mathbf{l}}},
   \end{equation}
   where $r= rank  M(n)$, $c_1, \ldots, c_r$ are positive numbers and $\mathbf{\lambda}_{\mathbf{l}} \in \mathbb{R}^d$, whenever $\mathbf{l} \in \mathfrak{I}$. By virtue of \cite[Theorem 7.8]{CuFi1}, $M(n+1)$ admits a unique (positive) flat extension $M(\infty) \equiv M(\beta)$, that is, $\beta^{(2n)}$ is a subsequence of some recursively generated moment sequence  $\beta \equiv \{\beta_\mathbf{i}\}_{ \mathbf{i} \in \mathbb{Z}_+^d}$.

 For  $k \in \{1, \ldots, m\}$, denote  $q_k (t_1, \ldots, t_d) = \sum\limits_{\mathbf{\alpha}} q_{k, \mathbf{\alpha}} t_1^{\alpha_1} \ldots t_d^{\alpha_d}$. From \eqref{rl2}, we have
  $\beta_{(i_1, \ldots, i_d)}= \sum\limits_{\mathbf{l} \in \mathfrak{I}} c_{\mathbf{l}} \lambda_{1, \mathbf{l}}^{i_1} \ldots \lambda_{d, \mathbf{l}}^{i_d}$; then
  \begin{equation}\label{rl3}
  \begin{aligned}
  (q_k * \beta)_{(i_1, \ldots, i_d)} &= \sum\limits_{\mathbf{\alpha}} q_{k, \mathbf{\alpha}} \sum\limits_{\mathbf{l} \in \mathfrak{I}} c_{\mathbf{l}} \lambda_{1, \mathbf{l}}^{i_1 +\alpha_1} \ldots \lambda_{d, \mathbf{l}}^{i_d +\alpha_d} \\
    &= \sum\limits_{\mathbf{l} \in \mathfrak{I}} c_{\mathbf{l}} q_k(\lambda_{1, \mathbf{l}}, \ldots, \lambda_{d, \mathbf{l}}) \lambda_{1, \mathbf{l}}^{i_1} \ldots \lambda_{d, \mathbf{l}}^{i_d}\\
    &= \sum\limits_{\mathbf{l} \in \mathfrak{I}} c_{\mathbf{l}} q_k(\mathbf{\lambda}_\mathbf{l}) \mathbf{\lambda}_\mathbf{l}^\mathbf{i}.
  \end{aligned}
  \end{equation}
  As
   \begin{equation*}
  \begin{aligned}
  c_{\mathbf{l}} q_k(\lambda_{1, \mathbf{l}}, \ldots, \lambda_{d, \mathbf{l}}) &= \int \mid L_{(\lambda_{1, \mathbf{l}} , \ldots, \lambda_{d, \mathbf{l}})} \mid^2 q_k d\mu \\
  &= L_{(\lambda_{1, \mathbf{l}} , \ldots, \lambda_{d, \mathbf{l}})}^T M_{q_k}(\beta) L_{(\lambda_{1, \mathbf{l}} , \ldots, \lambda_{d, \mathbf{l}})},
  \end{aligned}\end{equation*}
  and $M_{q_k}(\beta) \equiv M(q_k *\beta) \geq 0$, then $c_{\mathbf{l}} q_k(\lambda_{1, \mathbf{l}}, \ldots, \lambda_{d, \mathbf{l}}) \geq 0$. Thus,
   for every  $\mathbf{l}= (l_1, \ldots, l_d) \in \mathfrak{I}$, we obtain
  $q_k(\lambda_\mathbf{l}) = q_k(\lambda_{1, \mathbf{l}}, \ldots, \lambda_{d, \mathbf{l}}) \geq 0$ and this implies that  $supp \mu \subseteq K_{\mathcal{Q}}$, as desired.

  Since $\beta$ is recursively generated,  $M(\beta) \geq 0$ and $rank M(\beta) = rank M(n)$, we derive from Relation \eqref{rl2} and Theorem \ref{11} that $\mu = \sum\limits_{\mathbf{l} \in \mathfrak{I}} c_{\mathbf{l}} d\delta_{\lambda_{1, \mathbf{l}}} \ldots d\delta_{\lambda_{d, \mathbf{l}}}$ is the unique representing measure of $\beta$. Similarly, since for every $i =1, \ldots, m$,  $\{ (q_i*\beta)_\alpha \}_{\alpha \in \mathbb{Z}_+^d}$ is a recursively generated sequence, with $rank  M(q_i *\beta) =  rank  M_{q_i} (\infty) = rank  M_{q_i} (n +[\frac{\deg q_i +1}{2}])$ and $M_{q_i} (n +[\frac{\deg q_i +1}{2}]) \geq 0$. We get, by applying Theorem \ref{11} and Relation \eqref{rl3}, that  $\{(q_i * \beta)_\alpha  \}_{\alpha \in \mathbb{Z}_+^d}$ admits a unique representing measure
  \begin{equation}\label{mu2}
  \mu_i = \sum\limits_{\mathbf{l} \in \mathfrak{I}} c_{\mathbf{l}} q_i (\lambda_{1, \mathbf{l}}, \ldots, \lambda_{d, \mathbf{l}}) d\delta_{\lambda_{1, \mathbf{l}}} \ldots d\delta_{\lambda_{d, \mathbf{l}}},
   \end{equation}
   which is $M_{q_i} (n +[\frac{\deg q_i +1}{2}])$-atomic. Thus, from \eqref{mu1} and \eqref{mu2}, $\mu$ has precisely $rank M(n) - rank M_{q_i}(n + [\frac{\deg q_i +1}{2}])$ atoms in $\mathcal{Z}(q_i) := \{ t\in \mathbb{R}^d: q_i(t)=0 \}$, for every  $1\leq i\leq m$.\fin


\begin{thebibliography}{99}
\bibitem{BT}  C. Bayer and J. Teichmann, The proof of Tchakaloff's theorem Proc. Amer. Math. Soc., 134:10 (2006), 3035-3040.
\bibitem{BRZ} R. Ben Taher, M. Rachidi and H. Zerouali, Recursive subnormal completion and truncated moment problem, Bull. London Math. Soc. 33 (2001) 425-432.

\bibitem{CuFi1} R. E. Curto and L. A. Fialkow, solution of the truncated complex moment problem for flat data, Mem. Amer. Math. Soc. 119 (1996).
\bibitem{CuFi10}  R. E. Curto and L. A. Fialkow, Truncated $K$-moment problems in several variables, J. Operator Theory 54 (2005), 189- 226.
\bibitem{CuFi4} R.E. Curto and L.A. Fialkow. The truncated complex K-moment problem. Trans. Amer. Math. Soc. 352 (2000), 2825- 2855.
\bibitem{DRS} F. Dubeau, W. Motta, M. Rachidi, and O. Saeki, On weighted r-generalized Fibonacci sequences, The Fibonacci Quarterly, 35 (1997) 102-110.
\bibitem{Las1} J.B. Lasserre. Global optimization with polynomials and the problem of moments, SIAM J. Optim. 11, pp. 796-817.
\bibitem{Las2} J.B. J.B. Lasserre. Polynomials nonnegative on a grid and discrete optimization, Trans. Amer. Math. Soc. 354, pp. 631-649.
\bibitem{Las3} J.B. Lasserre. Moments, Positive Polynomials and Their Applications, Imperial College Press, 2009.
\bibitem{Lau}  J.B. Laurent, M. (2005). Revisiting two theorems of Curto and Fialkow on moment matrices, Proc. Amer. Math. Soc. 133, pp. 2965-2976.
\bibitem{Smu} J.L. Smul'jan, An operator Hellinger integral (Russian), Mat. Sb. 91 (1959), 381-430.


\end{thebibliography}
\end{document}